.
.
\font\script=eusm10.
\font\sets=msbm10.
\font\stampatello=cmcsc10.
\font\symbols=msam10.

\def\1{\hbox{\bf 1}}

\def\starsum{\mathop{{\sum}^{\ast}}}

\def\square{\hbox{\vrule\vbox{\hrule\phantom{s}\hrule}\vrule}}
\def\defineq{\buildrel{def}\over{=}}
\def\defin{\buildrel{def}\over{\Longleftrightarrow}}

\def\N{\hbox{\sets N}}

\def\R{\hbox{\sets R}}
\def\Z{\hbox{\sets Z}}
\def\Corr{\hbox{\script C}}
\def\EssBdd{\hbox{\symbols n}\,}

\par
\noindent
\centerline{\bf On the Selberg integral of the $k-$divisor function}
\centerline{\bf and the $2k-$th moment of the Riemann zeta function}
\medskip
\smallskip
\centerline{\stampatello Giovanni Coppola}
\bigskip
\par
\noindent {\bf 1. Introduction and statement of the result.}
\smallskip
\par
\noindent
We will link the $2k-$th moment of the {\stampatello Riemann $\zeta-$function} {\stampatello on the} (critical) {\stampatello line} ($\sigma={1\over 2}$), [Iv0]:
$$
I_k(T)\defineq \int_{0}^{T}|\zeta({\textstyle {1\over 2}} + it)|^{2k}dt
$$
\par
\noindent
(which we'll abbreviate with $I_k$, not to be confused with the similar $2k-$th moment {\stampatello off the line}, i.e.
$$
I_k(\sigma,T)\defineq \int_{0}^{T}|\zeta(\sigma+ it)|^{2k}dt, \qquad {1\over 2}<\sigma<1, 
$$
\par
\noindent
compare [Iv0]), with the {\stampatello Selberg integral} of the $k-$divisor function, $d_k(n)$ (having Dirichlet series $\zeta^k$) 
$$
J_k(x,h)\defineq \int_{hx^{\varepsilon}}^{x}\Big| \sum_{t<n\le t+h}d_k(n)-M_k(t,h)\Big|^2dt
$$
\par
\noindent
(compare [C]; abbreviate $J_k$, now on), where, say, $M_k(t,h)$ \thinspace is the \lq \lq {\stampatello expected value}\rq \rq \thinspace of the (inner) sum. 
\par
This gives $d_k$ over the (say) \lq \lq {\stampatello short interval}\rq \rq \thinspace $[t,t+h]$ (as $h=o(t)$ $\forall t\in [hx^{\varepsilon},x]$); here and in the sequel $\varepsilon>0$ will be arbitrarily small, not the same at each occurrence.
\par
Actually, Ivi\'c gave (in [Iv2], to appear on JTNB for JAXXV Proc.) a non-trivial bound for \thinspace $J_k(x,h)$ \thinspace when the {\stampatello width} of the {\stampatello s.i.} (abbrev. short interval), {\stampatello namely} \thinspace $\theta := {{\log h}\over {\log x}}$ \thinspace is greater than $\theta_k \defineq 2\sigma_k -1$ (with $\sigma_k$ {\stampatello Carlson's abscissa}, i.e. $\inf \{\sigma \in ]1/2,1[ \, : \, I_k(\sigma,T)\ll T\}$, here): 
$$
\theta > \theta_k \enspace \Rightarrow \enspace \exists \delta=\delta(k)>0 \thinspace : \thinspace J_k(x,h)\ll {{xh^2}\over {x^{\delta}}}
\leqno{\hbox{({\stampatello Ivi\'c}, {\rm JTNB})}}
$$
\par
\noindent
(with trivial bound : \enspace $J_k(x,h)\ll xh^2 (\log x)^c$, \enspace where $c=c(k)>0$, see the following). 
\par
This result clearly gives non-trivial bounds for \thinspace $J_k$, using $\zeta-$moments information (off the critical line). For example, \enspace $\theta_3={1\over 6}$, $\theta_4={1\over 4}$, $\theta_5={{11}\over {30}}$, $\ldots$ (from values of $\sigma_k$). 
\par
So, knowledge of the $\zeta$ implies knowledge of the $d_k$ in {\stampatello a.a.} (abbrev. almost all) the {\stampatello s.i.} (short intervals).
\smallskip
\par
\noindent
\centerline{(See: \thinspace $J_k$ \thinspace non trivial $\Rightarrow$ \enspace ${\displaystyle \sum_{t<n\le t+h} }d_k(n)\sim M_k(t,h)$, {\stampatello a.a.s.i.})}
\medskip
\par
\noindent
However, we can also go in the opposite direction: if we have some kind of non-trivial information about the $d_k$, we can improve our knowledge (at least, on the $2k-$th moments) of the Riemann $\zeta-$function. Actually, this idea is due to Ivi\'c, who linked $I_k$ to the \lq \lq {\stampatello (auto-)correlation}\rq \rq \thinspace of $d_k$ with \lq \lq {\stampatello shift-parameter}\rq \rq \thinspace $a$, i.e. 
$$
\Corr_k(a)\defineq \sum_{n\le x}d_k(n)d_k(n+a)\qquad a\in \N \enspace (\hbox{\rm here} \enspace x\in \N, x\to \infty)
$$
\par
\noindent
(the $\underline{\hbox{\stampatello shift}}$ is a positive integer: $\Corr_k(-a)$ is close enough to $\Corr_k(a)$ and $\Corr_k(0)$ is relatively easy to compute). 
\par
Here it comes into play the idea of Ivi\'c (see [Iv1] in Palanga 1996 Conference Proc.) of {\stampatello linking} the estimate of the {\stampatello $2k-$th moment}, $I_k(T)$, {\stampatello to a sum of correlations} $\Corr_k(a)$ performed {\stampatello over} $a$ ({\stampatello the shift}), up to (roughly, we avoid technicalities), say, $h:={x\over T}$ (the {\stampatello s.i.} comes in !)
\medskip
\par
In order to be more precise, we {\stampatello need to abbreviate} (with \thinspace $x,X$ \thinspace or even \thinspace $T$ \thinspace our \lq \lq {\stampatello main variables}\rq \rq, all independent \& $\to \infty$):
$$
A\EssBdd B \defin \forall \varepsilon>0 \enspace A\ll_{\varepsilon} x^{\varepsilon}B
$$
\par				% PAGE
\noindent
{\stampatello i.e.,} the {\stampatello modified Vinogradov notation} $\EssBdd$ allows us to {\stampatello ignore all the arbitrarily small powers}; also, we'll say that the {\stampatello arithmetic function} ({\stampatello a.f.}) $f:\N \rightarrow \R$ is {\stampatello essentially bounded}, write $f\EssBdd 1$, when $\forall \varepsilon>0$ $f(n)\ll_{\varepsilon} n^{\varepsilon}$ (as $n\to \infty$). For example, all the $d_k$ ($\forall k\in \N$) are ess.bd. :
$$
\forall k\in \N \quad d_k\EssBdd 1
$$
\par
\noindent
whence they contribute individually small powers (ignored); Shiu [S] estimates (see $J_k$ triv.est.quoted above), a kind of Brun-Titchmarsh for (suitable multiplicative {\stampatello a.f.}, like) $d_k$, let these give, on average over (all) {\stampatello s.i.}, powers of $\log$. By the way,
$$
L:=\log x \quad (\hbox{\rm or} \enspace L:=\log X, \enspace \hbox{\rm even} \enspace L:=\log N)
$$
\par
\noindent
is the abbreviation for the logarithm of our main variable.
\medskip
\par
We quote the formula (proved $\forall k\le 2$, see $\S3$) for $d_k$ correlations 
$$
\Corr_k(a)=x P_{2k-2}(\log x) + \Delta_k(x,a), \enspace \Delta_k(x,a)=o(x);
\leqno{(*)_k}
$$
\par
\noindent
here, the (conjectured, $\forall k>2$) {\stampatello main term} of $(\ast)_k$ is \enspace $xP_{2k-2}(\log x)\ll_k xL^{2k-2}\EssBdd x$ (since $P_{2k-2}$ is a polynomial of deg.$2k-2$, see the following). 
\par
Here, it seems that the first to propose explicitly this form for $(\ast)_k$ is Ivi\'c, who also gave explicitly the polynomial \enspace $P_{2k-2}$, that is essentially bounded (w.r.t. $x$). However, as we'll see in a moment, it depends, also, on the {\stampatello shift} $a>0$.
\medskip
\par
We'll give now, avoiding technicalities, Ivi\'c' s argument.
\medskip
\par
\noindent
After some work (expand the square \& mollify, take relevant ranges, $\ldots$) he gets that $I_k(T)$ is 
$$
I_k(T)=I''_k(T)+{\cal O}_{\varepsilon}(T^{\varepsilon}T)
$$
\par
\noindent
with 
$$
I''_k\defineq {1\over M}\sum_{a\le h}\sum_{M<n\le M'}d_k(n)d_k(n+a)\int_{T\over 2}^{2T}\phi(t)e^{ita/n}dt,
$$
\par
\noindent
where $M<M'\le 2M$, with $M\EssBdd T^{k/2}$, say $h\EssBdd M/T$, the smooth (i.e., $C^{\infty}$) test-function $\phi$ has support into $]T/2,2T[$, $\phi([3T/4,4T/3])\equiv 1$, and has good decay 
$$
\phi^{(R)}(t)\ll_R T^{-R}, \enspace \forall R\in \N.
$$ 
\smallskip
\par
Now on (see the reason in next section) we can ignore (in bounds for $I_k$) all terms which are\thinspace $\EssBdd T$.
\smallskip
\par
\noindent
We give an idea of the polynomial, $P_{2k-2}$, given by Ivi\'c, before to proceed. It's (see [Iv1] for details)
$$
P_{2k-2}(\log x)\defineq {1\over x}\int_{0}^{x}\sum_{q=1}^{\infty}{{c_q(a)}\over {q^2}}R_k^2(\log t)dt,
$$
\par
\noindent
with, say, 
$$
R_k(\log t)\defineq {{C_{-k}(q)}\over {(k-1)!}}\log^{k-1}t + {{C_{1-k}(q)}\over {(k-2)!}}\log^{k-2}t + \ldots + {{C_{-2}(q)}\over {1!}}\log t + C_{-1}(q)
$$
\par
\noindent
depending on $q$, but not on $a$ (this is vital); also, w.r.t. $x$, $R_k(\log t)\EssBdd 1$ and this is very useful ! We'll see in a moment that the shape of these $C_{j}(q)$ is important only in case $q=1$. By the way, here \thinspace $c_q(a)$ \thinspace is the {\stampatello Ramanujan sum}, defined as ($\starsum$ is over $q-$coprime $j$s)
$$
c_q(a)\defineq \starsum_{j(\!\!\bmod q)}e_q(ja)=\sum_{{d|q}\atop {d|a}}d\mu\Big({q\over d}\Big)
$$
\par				% PAGE
\noindent
Hence, say, $S(a)\defineq \max(0,h-|a|)$ (here $\hat{S}$ is Fej\'er's kernel) gives 
$$
\widehat{S}\Big( {j\over q}\Big)\defineq \sum_a S(a)e_q(ja)\ge 0 \Rightarrow \sum_a S(a)c_q(a)\ge 0
$$
\par
\noindent
(see the link with $J_k$ soon) and from the elementary, $\forall d\in \N$, 
$$
\sum_{a\atop {a\equiv 0(d)}}S(a)=h+2\sum_{b\le h/d}(h-db)={{h^2}\over d}+d\left\{ {h\over d}\right\}\left( 1-\left\{ {h\over d}\right\}\right), 
$$
\par
\noindent
we get (apply $c_q(a)$, above), writing $\1_\wp=1$ if $\wp$ holds, $=0$ else: 
$$
\sum_a S(a)c_q(a)=\1_{q=1}h^2+\sum_{d|q}d^2\mu\Big({q\over d}\Big)\left\{ {h\over d}\right\}\left( 1-\left\{ {h\over d}\right\}\right).
$$
\par
\noindent
(It is here evident $q=1$ greater importance.) Thus, (see Ivi\'c [Iv1] and compare [C]): 
$$
\sum_a S(a) x P_{2k-2}(\log x)=h^2 \int_{hx^{\varepsilon}}^{x}R^2_k(1,\log t)dt + \hbox{\stampatello tails},
\leqno{(1)}
$$
\par
\noindent
where we mean, by \lq \lq {\stampatello tails}\rq \rq, remainders which are $\EssBdd h^3$. Here, the part of \thinspace $R_k(\log t)$ \thinspace term with \thinspace $q=1$ \thinspace is, say, 
$$
R_k(1,\log t)\defineq {{C_{-k}(1)}\over {(k-1)!}}\log^{k-1}t + {{C_{1-k}(1)}\over {(k-2)!}}\log^{k-2}t + \ldots + {{C_{-2}(1)}\over {1!}}\log t + C_{-1}(1)
$$
\par
\noindent
and gives (see the above) the term $M_k(\log t)$ into the Selberg integral; as it should be, since (from an elementary version of Linnik's Dispersion method, compare [C] Lemmas), assuming $(\ast)_k$ with this $P_{2k-2}$, we get 
$$
J_k(x,h)\sim \sum_a S(a)\Corr_k(a)-h^2 \int_{hx^{\varepsilon}}^{x}M^2_k(\log t)\,dt 
\sim \sum_a S(a)\Delta_k(x,a),
\leqno{(2)}
$$
\par
\noindent
where $\sim$ means ignoring \lq \lq {\stampatello tails}\rq \rq \thinspace (see above) \& \lq \lq {\stampatello diagonals}\rq \rq, i.e. remainders $\EssBdd xh$. We remark that {\stampatello both these errors are negligible} (at least, for $k=3,4$, see $\S5$), since they both contribute $\EssBdd T$ to \thinspace $I_k(T)$. 
\par
Then, due to $I''_k$ expression, Ivi\'c [Iv1] made a hypothesis about (avoiding technicalities) sums of $\Delta_k(x,a)$ (remainders into $(\ast)_k$ above), performed over the shift $a$, say $G_k$, which implies the bound \thinspace $I_k(T)\EssBdd T$ (for the same $k>2$). {\stampatello Now on} $k>2$.
\smallskip
\par
Of course, he doesn't need $(\ast)_k$ to hold $\underline{\hbox{\stampatello individually}}$ $\forall a$ ($\le h$, here), but he observes that he's summing up, into $G_k$, $\underline{\hbox{\stampatello without}}$ the modulus over the remainder, $\Delta_k(x,a)$, so some $a-$cancellation can take place.
\smallskip
\par
So far, he passes from an asymptotic formula $(\ast)_k$ to an $a-$averaged form of it, which is easier to prove (however, yet nobody has done it !).
\bigskip
\par
Here, with applications in mind,
\smallskip
\par
\centerline{\stampatello we pass from a single average to a $\underline{\underline{\hbox{\stampatello double}}}$ average}
\medskip
\par
\noindent
Building on his expression for $I''_k$, it's possible to make a less stringent hypothesis, to {\stampatello have a more flexible procedure} for the remainders $\Delta_k(x,a)$.
\smallskip
\par
We use, also, our previous work on the Selberg integral of the {\stampatello a.f.} $f$ ({\stampatello essentially bounded \& real}), compare [C], in order to let the Selberg integral of $d_k$, i.e. $J_k(x,h)$, come into play. (It is a kind of \lq \lq {\stampatello double average}\rq \rq \thinspace of $\Delta_k(x,a)$.)
\smallskip
\par
Unfortunately, due to an exponential factor multiplying $d_k(n)d_k(n+a)$ into $\Corr_k(a)$ we can't get a link with $I_k(T)$ using only $J_k(x,h)$ ({\stampatello with} $h\EssBdd {x\over T}$, $x\EssBdd T^{k/2}$), but we need, also, to make an hypothesis on another double average of remainders $\Delta_k(x,a)$. We give our Theorem and the proof in $\S4$.
\smallskip
\par				% PAGE
\noindent {\stampatello Theorem.} {\it Let } $M<M'\le 2M$, $T^{1+\varepsilon}\le M\ll T^{k/2}$ {\it and } $H=M^{1+\varepsilon}/T$, {\it with double average } $\widetilde{G}_k=\widetilde{G}_k(M,T)$ {\it defined as} 
$$
\widetilde{G}_k\defineq \sup_{M\le x\le M',t\le H}\left( {1\over t}J_k(x,t)+{1\over t}\left| \sum_{h\le t}\sum_{h<a\le t}\Delta_k(x,a)\right|\right).
$$
\par
\noindent
{\it Then for } $k=3,4$ \thinspace {\it and any fixed } $\varepsilon>0$ {\it we have} 
$$
I_k(T)\EssBdd T\left( 1+\sup_{T\EssBdd M\ll T^{k/2}}\widetilde{G}_k(M,T)/M\right).
$$
\medskip
\par
\noindent
In next two sections we'll briefly mention some history of $I_k$ and the (related) additive divisor problems. Then, we'll prove our Theorem in the subsequent section and, finally, we'll give some remarks in the fifth.

\bigskip

\par
\noindent {\bf 2. A concise history of the Riemann-zeta moments (on the line)...}
\smallskip
\par
\noindent
We should keep in mind, here, that for fixed $k\in \N$ we seek 
$$
I_k(T)=\int_{0}^{T}|\zeta({\textstyle {1\over 2}} + it)|^{2k}dt \EssBdd T \enspace (\hbox{\stampatello \lq \lq $2k-$th\thinspace moment\thinspace pbm\rq \rq})
$$
\par
\noindent
(we call it \lq \lq {\stampatello on the line}\rq \rq, since $\sigma={1\over 2}$ is the critical line; \lq \lq {\stampatello off the line}\rq \rq \thinspace means with ${1\over 2}<\sigma<1$ or $\sigma=1$)
\par
\noindent
that (for $k>2$) is our aim; in fact, English school gave first $2$ cases: first Hardy \& Littlewood [H-L] in 1916 gave asymptotics for $k=1$ (not too difficult!) and then Ingham [In] in 1927 for $k=2$ (actually, for both $k=1,2$ he gave only $P_{2k-2}$ leading term, hence error $\log x$ better than main term); then, Heath-Brown in 1979 [HB] gave, for $k=2$ (solved $2-$add.div.pbm., i.e. the binary additive divisor pbm, see $\S3$, using Weil's bound for Kloosterman sums), $P_{6}$ (not explicitly!) plus error $E_2\EssBdd T^{7/8}$.
\par
Starting from '94 \& '95, a series of Ivi\'c \& Motohashi papers (applying $SL(\Z,2)$ considerations for the binary add.div.pbm.) gave $E_2\EssBdd T^{2/3}$ and, in mean-square, even $E_2\EssBdd \sqrt{T}$. (Here log-pows, not small pows!). Ivi\'c explicited $P_{2k-2}$ (when $k=2$). Like the binary add.div.pbm., this is not the whole story ! 
\par
The case $k=3$, again (recall $\Corr_3(a)$ {\stampatello pbm}), is unsolved.
\par
\noindent
The bound 
$$
I_3(T)\EssBdd T
$$
\par
\noindent
is called {\stampatello the \lq \lq sixth moment\rq \rq \thinspace pbm} (actually, this is the {\stampatello weak version}) \& has a {\stampatello link} (in Ivi\'c, Proc. Cardiff 1996 Symposium) {\stampatello with} the {\stampatello ternary additive divisor pbm}.
\medskip
\par
\noindent
Another interesting moment (Heath-Brown '79): $I_{6}\ll T^2L^c$.
\medskip
\par
One glimpse, to the {\stampatello high moments} (instead, for $k\le 2$, see [I-M]): {\stampatello predicted asymptotics} is \thinspace $I_k \sim C(k)TL^{k^2}$, $\forall k\ge 1$ (English school,again!) applying {\stampatello Random Matrix Theory} (Quantum Physics concepts inspiration!) in (2000s) seminal works of J. Keating \& N. Snaith (at Bristol); many others (Conrey, Ghosh to name two). The RMT$-\zeta$ link originated (Dyson-Montgomery coffee-break) in 1972. Not the end$\ldots$

\bigskip

\par
\noindent {\bf 3. ... and of (some) additive divisor problems.}
\smallskip
\par
\noindent
The problem of proving $(\ast)_k$ (at least {\stampatello fixed} $a>0$) is the $\underline{\hbox{\stampatello $k-$ary additive divisor problem}}$: trivial case $k=1$ ($\Corr_1(a)=x$ $\forall a \in \Z$) \& the \lq \lq {\stampatello binary additive divisor problem}\rq \rq, $k=2$, are the only solved pbms.
\par
Case $k=3$ is the $\underline{\hbox{\stampatello ternary additive divisor problem}}$ (sometimes called \lq \lq {\stampatello Linnik problem}\rq \rq): some time ago, Vinogradov \& Takhtadzhjan (see below $k=2$) announced its solution but with, as yet, unfilled holes in their (extremely technical !) \lq \lq proof\rq \rq. Their approach still suffers from our lack of information about $SL(\Z,3)$; while our (enough good) state of the art about, instead, $SL(\Z,2)$ (actually, through Kuznetsov Formula application, see [T-V]) allowed (starting from [M] approach) Ivi\'c, Motohashi and Jutila to solve satisfactorily, see esp. [I-M] (and the recent Meurman's [Me]), the binary additive divisor problem (different approaches work, with weaker remainders). We mention (still $k=2$), in passing, Kloosterman sums bounds (like Weil's) in the $\delta-$method of Duke-Friedlander-Iwaniec for \lq \lq determinantal equations\rq \rq (esp.,[DFI]). An even more general problem than this last has been solved by Ismoilov (see Math.Notes 1986).
\par				% PAGE
Thus, so far, no one has proved (for $k>2$), given $a\in \N$, 
$$
\Corr_k(a)=x P_{2k-2}(\log x) + \Delta_k(x,a),\enspace \Delta_k(x,a)=o(x),
$$
\par
\noindent
{\stampatello as \enspace $x\to \infty$ \enspace (the \enspace $k-$ary \enspace additive \enspace divisor \enspace problem),} not even for a single shift $a>0$ (already $k=2$ has delicate \lq \lq $a-${\stampatello uniformity}\rq \rq \thinspace issues: [I-M]).
{\stampatello (main term's $P_{2k-2}(\log x)$'s a $2k-2$ deg. $\log x$ poly)}

\bigskip

\par
\noindent {\bf 4. Proof of the Theorem.}
\smallskip
\par
\noindent
First of all, main terms in $(1)$, with $P_{2k-2}$, are treated like Ivi\'c does [Iv1]; he has, in partial summations, {\stampatello say}, 
$$
\sum_{a\le t}\Delta_k(x,a),\enspace \hbox{\stampatello which is}\enspace 
{1\over t}\sum_{h\le t}\sum_{a\le t}\Delta_k(x,a)
={1\over t}\sum_{h\le t}\sum_{a\le h}\Delta_k(x,a)+{1\over t}\sum_{h\le t}\sum_{h<a\le t}\Delta_k(x,a), 
$$
\par
\noindent
where the second (double) sum is in our $\widetilde{G_k}$; the former is the {\stampatello arithmetic mean}
$$
{1\over t}\sum_{h\le t}\sum_{a\le h}\Delta_k(x,a) 
$$ 
\par
\noindent
(a kind of average, something like $C^1$ process in Fourier series) and can be expressed as (exchanging sums)
$$
{1\over t}\sum_{a\le t}(t-a+1)\Delta_k(x,a)={1\over t}\sum_{a\le t}(t-a)\Delta_k(x,a)+{1\over t}\sum_{a\le t}\Delta_k(x,a)
$$ 
\par
\noindent
which reduces to (using \enspace $\Delta_k(x,0)\EssBdd x$ \enspace for {\stampatello diagonals} and \enspace $\Delta_k(x,-a)=\Delta_k(x,a)+{\cal O}_{\varepsilon}(x^{\varepsilon}a)$ \enspace for {\stampatello tails}) 
$$
\sim {1\over t}\sum_{a\le t}(t-a)\Delta_k(x,a)\sim {1\over {2t}}\sum_{0\le |a|\le t}(t-|a|)\Delta_k(x,a)
$$
\par
\noindent
($+$ {\stampatello diagonals} \& {\stampatello tails}); and, since \thinspace $S(a)=\max(t-|a|,0)$, $\forall 0\le |a|\le t$ (apply $(2)$ \& compare [C]) $\Rightarrow $
$$
\sum_{0\le |a|\le t}(t-|a|)\Delta_k(x,a)\sim J_k(x,t), 
$$
\par
\noindent
we get the desired bound with Selberg integral (and double average).$\enspace \square$

\bigskip

\par
\noindent {\bf 5. Remarks.}
\smallskip
\par
\noindent
We remark that, in spite of the fact that Ivi\'c 's Theorem [Iv1] holds $\forall k>2$ (integer), we have some trouble in handling Selberg's integral {\stampatello tails}, since they contribute to \thinspace $\widetilde{G}_k$ \thinspace as (in the $\sup$ above)
$$
\EssBdd {1\over t}t^3 \EssBdd H^{2} \EssBdd {{M^2}\over {T^2}} 
$$
\par
\noindent
which gives to \thinspace $I_k(T)$ \thinspace a contribute (other $\sup$ above) 
$$
\EssBdd T\left( {{M^2}\over {T^2}}M^{-1}\right) \EssBdd {M\over T} \EssBdd T^{k/2-1} 
$$
\par
\noindent
that is \thinspace $\EssBdd T$ \thinspace {\stampatello only when} \thinspace $k/2\le 2$, {\stampatello i.e.} $k\le 4$ {\stampatello here}. 
\par
We trust the possibility to have a link as above not only for the sixth \& the eighth moment, but the {\stampatello tails} arise naturally when applying the Linnik method and even a more careful analysis will almost surely not eliminate them ! While they are negligible for the Selberg integral, not so for the present application !
\bigskip
\par
\noindent
We remark, in passing, that here the \lq \lq additional\rq \rq \thinspace double average can't be dispensed with.

\bigskip

\par				% PAGE
\centerline{\bf References}
\medskip
\item{\bf [C]} \thinspace Coppola, G.\thinspace - \thinspace {\sl On the Correlations, Selberg integral and symmetry of sieve functions in short intervals} \thinspace - \thinspace http://arxiv.org/abs/0709.3648v3
\smallskip
\item{\bf [D]} \thinspace Davenport, H.\thinspace - \thinspace {\sl Multiplicative Number Theory} \thinspace - \thinspace Third Edition, GTM 74, Springer, New York, 2000. $\underline{{\tt MR\enspace 2001f\!:\!11001}}$
\smallskip
\item{\bf [DFI]} Duke, W.,\thinspace Friedlander, J.\thinspace and \thinspace Iwaniec, H. \thinspace - \thinspace {\sl A quadratic divisor problem} \thinspace - \thinspace Inv. Math. {\bf 115} (1994), 209--217. 
\smallskip
\item{\bf [HB]} Heath-Brown, D. R.\thinspace - \thinspace {\sl The fourth power moment of the Riemann zeta function.} \thinspace - \thinspace Proc. London Math. Soc. (3) {\bf 38} (1979), 385--422. $\underline{{\tt MR81f\!:\!10052}}$
\smallskip
\item{\bf [H-L]} Hardy, G. H.\thinspace and \thinspace Littlewood, J. E. \thinspace - \thinspace {\sl Contributions to the theory of the riemann zeta-function and the theory of the distribution of primes} \thinspace - \thinspace Acta Math. {\bf 41} (1916), no. 1, 119--196. $\underline{{\tt MR1555148}}$
\smallskip
\item{\bf [In]} Ingham, A. E.\thinspace - \thinspace {\sl Mean-value theorems in the theorie of the Riemann Zeta-function} \thinspace - \thinspace Proc. London Math. Soc. (2) {\bf 27} (1927), 273--300. 
\smallskip
\item{\bf [Iv0]} \thinspace Ivi\'c, A.\thinspace - \thinspace {\sl The Riemann Zeta Function} \thinspace - \thinspace John Wiley \& Sons, New York, 1985. (2nd ed., Dover, Mineola, N.Y. 2003).
\smallskip
\item{\bf [Iv1]} \thinspace Ivi\'c, A.\thinspace - \thinspace {\sl The general additive divisor problem and moments of the zeta-function} \thinspace - \thinspace - New trends in probability and statistics, Vol. 4 (Palanga, 1996), 69--89, VSP, Utrecht, 1997. $\underline{{\tt MR\enspace 99i\!:\!11089}}$
\smallskip
\item{\bf [Iv2]} \thinspace Ivi\'c, A.\thinspace - \thinspace {\sl On the mean square of the divisor function in short intervals} \thinspace - \thinspace available online at the address http://arxiv.org/abs/0708.1601v2 \thinspace - \thinspace to appear on Journal de Th\'eorie des Nombres de Bordeaux.
\smallskip
\item{\bf [I-M]} \thinspace Ivi\'c, A.\thinspace and \thinspace Motohashi, Y.\thinspace - \thinspace {\sl On some estimates involving the binary additive divisor problem} \thinspace - \thinspace Quart. J. Math. Oxford Ser. (2) {\bf 46} (1995), no. 184, 471--483. - Proc. Cambridge (1995). $\underline{{\tt MR \enspace 96k\!:\!11117}}$
\smallskip
\item{\bf [Me]} \thinspace Meurman, T.\thinspace - \thinspace {\sl On the binary additive divisor problem} \thinspace - \thinspace Number theory (Turku, 1999), 223--246, de Gruyter, Berlin, 2001. $\underline{{\tt MR \enspace 2002f\!:\!11129}}$
\smallskip
\item{\bf [M]} \thinspace Motohashi, Y.\thinspace - \thinspace {\sl The binary additive divisor problem} Ann. Sci. \'Ecole Norm. Sup. (4) {\bf 27} (1994), no. 5, 529--572. $\underline{{\tt MR \enspace 95i\!:\!11104}}$
\smallskip
\item{\bf [S]} \thinspace Shiu, P. \thinspace - \thinspace {\sl A Brun-Titchmarsh theorem for multiplicative functions}  - J. Reine Angew. Math. {\bf 313} (1980), 161--170. $\underline{{\tt MR \enspace 81h\!:\!10065}}$
\smallskip
\item{\bf [T-V]} \thinspace Takhtadzhjan, L.A.\thinspace and \thinspace Vinogradov, A.I.\thinspace - \thinspace {\sl The zeta-function of the additive divisor problem and the spectral expansion of the automorphic Laplacian} (Russian) \thinspace - \thinspace Zap. Nauchn. Semin. Leningr. Otd. Math. Inst. Steklova {\bf 134} (1984), 84--116.
\medskip
\leftline{\tt Dr.Giovanni Coppola}
\leftline{\tt DIIMA - Universit\`a degli Studi di Salerno}
\leftline{\tt 84084 Fisciano (SA) - ITALY}
\leftline{\tt e-mail : giocop@interfree.it}

\bye